\documentclass[12pt]{article}
\usepackage{amsmath, amsthm}
\usepackage{amscd}
\usepackage{amssymb}
\usepackage{array}
\usepackage{longtable}
\usepackage{multirow}
\usepackage{color}
\newtheorem{thm}{Theorem}[section]
\newtheorem{theorem}[thm]{Theorem}

\newtheorem{lemma}[thm]{Lemma}
\newtheorem{proposition}[thm]{Proposition}

\theoremstyle{definition}
\newtheorem{definition}[thm]{Definition}

\begin{document}

\newcommand{\be}{\begin{equation}}
\newcommand{\ee}{\end{equation}}
\newcommand{\bea}{\begin{eqnarray}}
\newcommand{\eea}{\end{eqnarray}}
\newcommand{\bean}{\begin{eqnarray*}}
\newcommand{\eean}{\end{eqnarray*}}

\newcommand{\id}{\relax{\rm 1\kern-.28em 1}}
\newcommand{\R}{\mathbb{R}}
\newcommand{\C}{\mathbb{C}}
\newcommand{\Z}{\mathbb{Z}}
\newcommand{\g}{\mathfrak{G}}
\newcommand{\e}{\epsilon}

\newcommand{\hs}{\hfill\square}
\newcommand{\hbs}{\hfill\blacksquare}

\newcommand{\bp}{\mathbf{p}}
\newcommand{\bmax}{\mathbf{m}}
\newcommand{\bT}{\mathbf{T}}
\newcommand{\bU}{\mathbf{U}}
\newcommand{\bP}{\mathbf{P}}
\newcommand{\bA}{\mathbf{A}}
\newcommand{\bm}{\mathbf{m}}
\newcommand{\bIP}{\mathbf{I_P}}

\newcommand{\cA}{\mathcal{A}}
\newcommand{\cB}{\mathcal{B}}
\newcommand{\cC}{\mathcal{C}}
\newcommand{\cI}{\mathcal{I}}
\newcommand{\cO}{\mathcal{O}}
\newcommand{\cG}{\mathcal{G}}
\newcommand{\cJ}{\mathcal{J}}
\newcommand{\cF}{\mathcal{F}}
\newcommand{\cP}{\mathcal{P}}
\newcommand{\ep}{\mathcal{E}}
\newcommand{\E}{\mathcal{E}}
\newcommand{\cH}{\mathcal{O}}
\newcommand{\cPO}{\mathcal{PO}}
\newcommand{\cl}{\ell}
\newcommand{\cFG}{\mathcal{F}_{\mathrm{G}}}
\newcommand{\cHG}{\mathcal{H}_{\mathrm{G}}}
\newcommand{\Gal}{G_{\mathrm{al}}}
\newcommand{\cQ}{G_{\mathcal{Q}}}

\newcommand{\ri}{\mathrm{i}}
\newcommand{\re}{\mathrm{e}}
\newcommand{\rd}{\mathrm{d}}
\newcommand{\rid}{\mathrm{id}}

\newcommand{\rGL}{\mathrm{GL}}
\newcommand{\rSU}{\mathrm{SU}}
\newcommand{\rU}{\mathrm{U}}
\newcommand{\rSL}{\mathrm{SL}}
\newcommand{\rSO}{\mathrm{SO}}
\newcommand{\rOSp}{\mathrm{OSp}}
\newcommand{\rSpin}{\mathrm{Spin}}
\newcommand{\rsl}{\mathrm{sl}}
\newcommand{\rsu}{\mathrm{su}}
\newcommand{\ru}{\mathrm{u}}
\newcommand{\rM}{\mathrm{M}}
\newcommand{\rdiag}{\mathrm{diag}}
\newcommand{\rP}{\mathrm{P}}

\newcommand{\M}{\mathrm{M}}
\newcommand{\End}{\mathrm{End}}
\newcommand{\Hom}{\mathrm{Hom}}
\newcommand{\diag}{\mathrm{diag}}
\newcommand{\rspan}{\mathrm{span}}
\newcommand{\rank}{\mathrm{rank}}
\newcommand{\ber}{\mathrm{Ber}}

\newcommand{\fsl}{\mathfrak{sl}}
\newcommand{\fg}{\mathfrak{g}}
\newcommand{\fh}{\mathfrak{h}}
\newcommand{\ff}{\mathfrak{f}}
\newcommand{\fgl}{\mathfrak{gl}}
\newcommand{\fosp}{\mathfrak{osp}}
\newcommand{\fm}{\mathfrak{m}}
\newcommand{\fp}{\mathfrak{p}}
\newcommand{\fsu}{\mathfrak{su}}
\newcommand{\fu}{\mathfrak{u}}
\newcommand{\fn}{\mathfrak{n}}
\newcommand{\fb}{\mathfrak{b}}

\newcommand{\str}{\mathrm{str}}
\newcommand{\Sym}{\mathrm{Sym}}
\newcommand{\tr}{\mathrm{tr}}
\newcommand{\defi}{\mathrm{def}}
\newcommand{\Ber}{\mathrm{Ber}}
\newcommand{\spec}{\mathrm{Spec}}
\newcommand{\sschemes}{\mathrm{(sschemes)}}
\newcommand{\sschemeaff}{\mathrm{ {( {sschemes}_{\mathrm{aff}} )} }}
\newcommand{\rings}{\mathrm{(rings)}}
\newcommand{\Top}{\mathrm{Top}}
\newcommand{\sarf}{ \mathrm{ {( {salg}_{rf} )} }}
\newcommand{\arf}{\mathrm{ {( {alg}_{rf} )} }}
\newcommand{\odd}{\mathrm{odd}}
\newcommand{\alg}{\mathrm{(alg)}}
\newcommand{\sa}{\mathrm{(salg)}}
\newcommand{\sets}{\mathrm{(sets)}}
\newcommand{\SA}{\mathrm{(salg)}}
\newcommand{\salg}{\mathrm{(salg)}}
\newcommand{\varaff}{ \mathrm{ {( {var}_{\mathrm{aff}} )} } }
\newcommand{\svaraff}{\mathrm{ {( {svar}_{\mathrm{aff}} )}  }}
\newcommand{\ad}{\mathrm{ad}}
\newcommand{\Ad}{\mathrm{Ad}}
\newcommand{\pol}{\mathrm{Pol}}
\newcommand{\Lie}{\mathrm{Lie}}
\newcommand{\Proj}{\mathrm{Proj}}
\newcommand{\rGr}{\mathrm{Gr}}
\newcommand{\rFl}{\mathrm{Fl}}
\newcommand{\rPol}{\mathrm{Pol}}
\newcommand{\coinv}{\mathrm{coinv}}

\newcommand{\uspec}{\underline{\mathrm{Spec}}}
\newcommand{\uproj}{\mathrm{\underline{Proj}}}

\newcommand{\sym}{\cong}
\newcommand{\al}{\alpha}
\newcommand{\lam}{\lambda}
\newcommand{\de}{\delta}
\newcommand{\D}{\Delta}
\newcommand{\s}{\sigma}
\newcommand{\lra}{\longrightarrow}
\newcommand{\ra}{\rightarrow}


\newcommand{\beq}{\begin{equation}}
\newcommand{\eeq}{\end{equation}}
\newcommand{\Span}{\mathrm{span}}
\newcommand{\bM}{\mathbb{M}}



\bigskip



\bigskip

 \centerline{\LARGE \bf  $N=2$ quantum chiral superfields  }

    \smallskip

    \centerline{\LARGE \bf and quantum super bundles}

\vskip 1cm

\centerline{ R. Fioresi \footnote{Investigation supported by the
University of Bologna, funds for selected research topics.}
}
\smallskip
\centerline{\it Fabit, Universit\`a di
Bologna}

{\centerline{\it via San Donato, 15. 40126 Bologna, Italy}}

\centerline{\it and INFN,  Sezione di
    Bologna, Italy}

\centerline{{\footnotesize e-mail: rita.fioresi@unibo.it}}

\bigskip
\centerline{ M. A. Lled\'o }

\smallskip

 \centerline{\it  Departament de F\'{\i}sica Te\`{o}rica,
Universitat de Val\`{e}ncia and }

\centerline{\it IFIC (CSIC-UVEG)}
 \centerline{\it C/Dr.
Moliner, 50, E-46100 Burjassot (Val\`{e}ncia), Spain.}
 \centerline{{\footnotesize e-mail: maria.Lledo@ific.uv.es}}

\bigskip
\centerline{ J. Razzaq }

\smallskip
\centerline{\it Dipartimento di Matematica, Universit\`a di
Bologna,}
 \centerline{\it Piazza di Porta San Donato, 5. 40126 Bologna, Italy}

\centerline{{\footnotesize e-mail: junaid.razzaq2@unibo.it}}

\vskip 1cm

\begin{abstract}

 We give the {super}algebra of $N=2$ chiral (and antichiral) quantum superfields realized as a subalgebra of the quantum supergroup $\rSL_q(4|2)$. The multiplication law in the quantum supergroup induces a coaction on the set of chiral superfields. We also realize the quantum deformation of
the chiral Minkowski superspace as a quantum principal bundle.

\end{abstract}
\vfill\eject


\section{Introduction}\label{intro-sec}

It is well known that the $N=1$ superconformal superspace, in its complexified version \cite{ma1,va},  is the superflag  $Fl(2|0,2|1,4|1)$, on which the conformal supergroup $\rSL(4|1)$ acts naturally. The space $\C^{4|1}$,
{underlying}
the defining representation of
$\rSL(4|1)$, is the space of supertwistors \cite{pe,fe}.

Dealing with the complexified version has the advantage of seeing this structure, while the {conditions for the} real form can be imposed later on \cite{va}.  It is also well known, and differently from the non super case,
that not all the superflags are projective superspaces (take for example the super Grassmannian $Gr(1|1, 2|2)$ \cite {ma1}) and indeed the projective cases are rare among these superspaces, {though a new approach to this question was taken
  in \cite{sv}.}
For the super Grassmannians only the extreme cases
$Gr(p|0, m|n)$ or  $Gr(p|n, m|n)$ are superprojective and are both embedded into the projective superspace {$\mathbb{P}^{M|N}$ for suitable $M$ and $N$
see \cite{fl}}. These super Grassmannians are dual to each other and are the antichiral and chiral superspaces respectively.

 The superflag $Fl(2|0,2|1,4|1)$ can be embedded in the product
 $$Fl(2|0,2|1,4|1)\subset Gr(2|0,4|1)\times  Gr(2|1,4|1)\,,$$ and using the super Segre embedding \cite{flln} the superflag is embedded into the projective superspace $\mathbb{P}^{80|64}$ \cite{cfl,fl}.

 For $N=2$ we can reproduce the same situation with
 $$Fl(2|0,2|2,4|2)\subset Gr(2|0,4|2)\times  Gr(2|2,4|2)\,,$$ but this superflag is too big. The scalar superfields associated to it have too many field components to be useful in  the formulation of  supersymmetric field theories. Still, the antichiral $Gr(2|0,4|2)$ and chiral  $Gr(2|2,4|2)$ superspaces do have physical applications so it is useful to study them. They are both embedded in $\mathbb{P}^{8|8}$.

 There is a third super Grassmannian, $Gr(2|1,4|2)$, which is not projective but that has physical applications. It is the harmonic superspace of \cite{gikos,gios}. This example, together with more general ones,  were studied from this point of view in the series of papers \cite{hl,hh,hah,hho,heho,he,haho}.

 Here we will consider only the (anti)chiral superspace (also considered in \cite{me}). Our aim is to quantize it by substituting the supergroup $\rSL(4|2)$ by the  quantum group $\rSL_q(4|2)$ (in the sense of Manin \cite{ma2}) and trying to define appropriately the quantum super Grassmannian as an homogeneous superspace. This was done for $N=1$ in \cite{cfl,fl}. As we will see, the $N=2$ case has its own peculiarities.

 This program could, in principle be proposed for general homogeneous superspaces, not necessarily superprojective. But the projectivity gives us an advantage: the algebra associated to the projective embedding (super Pl\"{u}cker embedding) can be seen both, as a quotient algebra of the projective superspace 
$\mathbb{P}^{8|8}$ modulo some homogeneous polynomial relations (super Pl\"{u}cker relations) as well as a {graded} subring of the superring  $\C[\rSL(4|2)]$,
 {encoding its projective embedding (see also \cite{fi1, cfg, afl}).}
 We will see this in detail in Section \ref{classical-sec}.
  One can then define a quantum super Grassmannian  as a certain subalgebra of the super Hopf algebra $\rSL_q(4|2)$. If done correctly, the subalgebra must represent a quantum homogeneous superspace for  $\rSL_q(4|2)$, that is, the coproduct in  $\rSL_q(4|2)$ induces a coaction on the quantum super Grassmannian.

The chiral Minkowski $N=2$ superspace $\bM$ emerges naturally in
this context as the big cell in the Grassmannian  $Gr(2|0,4|2)$.
The $N=1$ case was extensively studied in \cite{fl}, Chapter 4. However, as
remarked above, the $N=2$ SUSY has its own peculiarities, which make
the theory richer. We view the big cell in $Gr(2|0,4|2)$ as
the subsupermanifold containing certain $2|0$ subspaces and we realize it
as the set $S$ of pairs of vectors
in $\C^{4|2}$ modulo the natural right $\rGL(2)$ action, which accounts
for basis change.
Hence, we construct $\bM$ as the quotient of $S$ modulo
the ordinary general linear group $\rGL(2)$.
The quantization of $\bM$ is obtained, as expected,
as the subsuperring of a localization of $\rSL_q(4|2)$,
generated by the quantum coinvariants with respect to the
coaction of quantum $\rGL_q(2)$ (see \cite{fl}, Chapter 4 for
the $N=1$ case). The presentation
of this quantum superring via generators and relations, makes an essential
use of the commutation relations among the quantum  determinants
appearing
in the definition of the quantum $Gr(2|0,4|2)$ and the
Pl\"{u}cker relations. Moreover, the quantum Minkowski  space, $\bM_q$, is isomorphic to
the quantum Manin superalgebra, that is, the quantum super bialgebra
of matrices, as described in \cite{ma2}. This fact is highly non
obvious, it depends on the quite involved
commutation relations of quantum determinants
and it shows how this framework is natural and suitable for
more exploration, as we detail below.

\medskip
The chiral  Minkowski $N=2$ superspace, being
a quotient, appears then naturally also
as a principal
bundle for the action of $\rGL(2)$. There is an extensive literature
regarding the quantization of principal bundles (see \cite{afl, hajac, hajac1,
hajac2, buachalla} and references therein). In particular the notion of
Hopf-Galois extension \cite{mont}
appears to be  the right one to formulate, in the
affine setting, the 
theory of principal bundles
to obtain their quantum deformations.

We hence proceed to define Hopf-Galois extensions in the
SUSY framework and prove that the chiral Minkowski $N=2$
superspace $\bM$ is {the base for
a principal bundle $S$} for the supergroup $\rGL(2)$,
by realizing it as a trivial Hopf-Galois extension
(see also \cite{kessler}  for a more geometric, yet equivalent, view
on super principal bundles).
Next, we construct a quantum deformation $\bM_q$ of $\bM$, by taking
advantage of our previous realization and show
that $\bM_q$ is the quantum space, {base for the 
quantum principal bundle $S_q$}, for $\rGL_q(2)$.

\medskip
We plan to explore, in a forthcoming paper, the construction
of covariant differential calculi on the quantum chiral Minkowski $N=2$
superspace and then proceed towards the realization of
a theory in a curved background.

\medskip

The paper is organized as follows.

In Section \ref{plucker-sec} we describe the super Pl\"{u}cker embedding of the super Grassmannian and its presentation in terms of generators and relations.

In Section \ref{classical-sec} we give the classical super  Grassmannian as a subalgebra of the coordinate superalgebra of $\rSL(4|2)$.

In Section \ref{bigcell-sec} we briefly describe the big cell of the super Grassmannian, the $N=2$, $D=4$ Minkowski superspace.

In Section \ref{quantumgrassmannian-sec} we pass to define the super Grassmannian  as a  subsuperalgebra of $\rSL_q(4|2)$, computing the commutation relations of the generators and the quantum super Pl\"{u}cker relations that they satisfy.

 In Section \ref{homogenous-sec} we give the coaction of $\rSL_q(4|2)$ on the quantum super Grassmannian defined in the previous section, proving that it is a quantum homogeneous superspace.

 In the last section, Section \ref{quantumbundle-sec},  we construct the $N=2$ ordinary chiral Minkowski superspace $\bM$ and its quantization $\bM_q$, realized both first as
homogeneous spaces for the action of the ordinary (quantum) general
linear group in dimension $2$, then as {bases for} (quantum)
principal bundles for  $\rGL(2)$ and  $\rGL_q(2)$ respectively.

 \section{The super Pl\"{u}cker embedding}\label{plucker-sec}

We are going to give the embedding of $Gr(2|0,4|2)$ in the projective superspace $\mathbb{P}^{8|8}$. 
Let $E=\bigwedge^2\C^{4|2}$  and $\{e_1,\dots,e_4,\epsilon_5,\epsilon_6\}$ an homogeneous basis for $\C^{4|2}$, we then have a basis for $E$ as
\begin{align*}
&e_i\wedge e_j \quad 1\leq i<j\leq 4, \quad \epsilon_5\wedge\epsilon _5,\quad \epsilon_6\wedge\epsilon _6,\quad \epsilon_5\wedge\epsilon _6,&\hbox{(even)},\\
&e_k\wedge \epsilon _5,\quad e_k\wedge\epsilon_6\quad 1\leq k\leq 4, &\hbox{(odd)}\,.
\end{align*}
So $E\simeq \C^{9|8}$ and $\mathbb{P}(E)\simeq\mathbb{P}^{8|8}$. An element of $E$ is given as
$$Q=q+\lambda_5\wedge\epsilon_5+\lambda_6\wedge\epsilon _6+a_{55}\epsilon_5\wedge\epsilon_5+a_{66}\epsilon_6\wedge\epsilon_6+ a_{56}\epsilon_5\wedge\epsilon_6\,,$$
with
$$q=q_{ij}e_i\wedge e_j,\quad \lambda_m=\lambda_{im}e_i,\qquad i,j=1,\dots, 4,\quad m=5,6\,.$$
The element $q$ is decomposable if $q=a\wedge b$, where
$$a=r+\xi_5\epsilon_5+\xi_6\epsilon_6,\qquad b=s+\eta_5\epsilon_5+\eta_6\epsilon_6,$$ with
$r=r_ie_i$, $s=s_ie_i$.

One obtains the following relations
\begin{align}
&q=r\wedge s&& &&\,\nonumber\\
&\lambda_5=\xi_5s-\eta_5 r,\qquad &&\lambda_6=\xi_6s-\eta_6 r,&&\nonumber\\
&a_{55}=\xi_5\eta_5,
\qquad &&a_{66}=\xi_6\eta_6,
\qquad &&a_{56}=\xi_5\eta_6+
\xi_6\eta_5\,,\label{wedgeproduct}
\end{align}
which imply
\begin{align}
&q\wedge q=0,\nonumber\\
&q\wedge\lambda_5=0,\qquad &&q\wedge\lambda_6=0,&&\nonumber\\
&\lambda_5\wedge\lambda_5=-2a_{55}q,\qquad &&\lambda_6\wedge\lambda_6=-2a_{66}q,\qquad &&\lambda_5\wedge\lambda_6=-a_{56}q,\nonumber\\
&\lambda_5 a_{55}=0,\qquad &&\lambda_6 a_{66}=0,\qquad&&\nonumber\\
&\lambda_5 a_{66}=-\lambda_6a_{56},\qquad &&\lambda_6 a_{55}=-\lambda_5a_{56},\qquad&&\nonumber
\\&a_{55}^2=0,\qquad &&a_{66}^2=0,\qquad && a_{56}a_{56}=-2a_{55}a_{66},\nonumber
\\&a_{55}a_{56}=0,\qquad &&a_{66}a_{56}=0\,.\qquad && \label{superplucker}
\end{align}

Relations (\ref{superplucker}) are  the super Pl\"{u}cker relations.
We can write them in coordinates in the following way (always $1\leq i<j<k\leq 4$ and $5\leq n\leq 6$):

\begin{align}
&q_{12}q_{34}- q_{13}q_{24}+q_{14}q_{23}=0,&&\hbox{(Pl\"{u}cker relation)}
\nonumber\\
&q_{ij}\lambda_{kn}-q_{ik}\lambda_{jn}+q_{jk}\lambda_{in}=0,
\nonumber\\
&\lambda_{in}\lambda_{jn}=a_{nn}q_{ij}, &&\lambda_{i5}\lambda_{j6}+\lambda_{i6}\lambda_{j5}=a_{56}q_{ij},
\nonumber\\
&\lambda_{in}a_{nn}=0,&&\lambda_{i5}a_{66}=-\lambda_{i6}a_{56}&&\lambda_{i6}a_{55}=-\lambda_{i5}a_{56}
\nonumber\\&a_{nn}^2=0&&a_{55}a_{56}=0,&&a_{66}a_{56}=0\nonumber\\
&a_{56}a_{56}=-2a_{55}a_{66}\,. &&\label{superpluckercoor}
\end{align}
For $N=1$ these relations were given in \cite{cfl}, were the relation $a_{55}^2=0$ was missing but implicitly assumed.  As we can see, for $N=2$  extra relations appear. The super Pl\"{u}cker relations for arbitrary $N$ are given in \cite{sv}, page 17. They coincide with ours  by a change in the notation
due to the appearance of a sign, because of a different convention
  on row/column vectors, hence the consequent change of sign
  of the supertranspose.

We will denote as $\cI_P$ the ideal generated by them in the affine superspace 
{$\mathbb{A}^{9|8}$} 
(with generators $q_{ij}, a_{nm}, \lambda_{kn}$).  They are homogeneous quadratic equations, so they are defined in the projective space $\mathbb{P}^{8|8}$.

Let us denote $\rGr= Gr(2|0,4|2)$ and consider the super Pl\"{u}cker map
$$\begin{CD}\rGr @>>>\mathbb{P}^{8|8}\\
\rspan\{a, b\}@>>>[a\wedge b]\,,\end{CD}$$

 We have the following

\begin{proposition}\label{superring-prop}The superring associated to the image of $\rGr$ under the super Pl\"{u}cker embedding is
$$\C[\rGr]\cong\C[q_{ij}, a_{nm}, \lambda_{kn}]/\cI_P\,,$$
that is, the relations in $\cI_P$ are all the relations satisfied by the generators $q_{ij}, a_{nm}, \lambda_{kn}$. Then $Gr(2|0,4|2)$, is a projective supervariety.\end{proposition}

{\it Proof.} This is  proven for arbitrary $N$ (and further generalizations) in \cite{sv}, Theorem 6, {denoted as the ``algebraic case''}.

\hfill$\blacksquare$

\section{ The classical picture}\label{classical-sec}

As stated in the introduction, one can see $\C[\rGr]$ as a subalgebra of
$\C[\rSL(4|2)]$. Let us display the generators of this algebra in matrix form
\be\begin{pmatrix}
g_{11}&g_{12}&g_{13}&g_{14}&\vline&\gamma_{15}&\gamma_{16}\\
g_{21}&g_{22}&g_{23}&g_{24}&\vline&\gamma_{25}&\gamma_{26}\\
g_{31}&g_{32}&g_{33}&g_{34}&\vline&\gamma_{35}&\gamma_{36}\\
g_{41}&g_{42}&g_{43}&g_{44}&\vline&\gamma_{45}&\gamma_{46}\\
\hline
\gamma_{51}&\gamma_{52}&\gamma_{53}&\gamma_{54}&\vline&g_{55}&g_{56}\\
\gamma_{61}&\gamma_{62}&\gamma_{63}&\gamma_{64}&\vline&g_{65}&g_{66}
\end{pmatrix}\,\label{groupmatrix}\ee then

$$\C[\rSL(4|2)]=\C[g_{ij},g_{mn}, \gamma_{im}, \gamma_{nj}]/(\mathrm{Ber}-1)\,,
$$
where $\mathrm{Ber}$ is the Berezinian of the matrix and $1\leq i,j\leq 4$ and  $5\leq m,n\leq 6$.

\begin{proposition}The superring\label{generators-prop}
$$\C[\rGr]\cong\C[q_{ij}, a_{nm}, \lambda_{kn}]/\cI_P\,,$$
is generated, as a subring of $\C[\rSL(4|2)]$ by the elements
\begin{align*}
&y_{ij}=g_{i1}g_{j2}-g_{i2}g_{j1},&&\eta_{kn}=g_{k1}\gamma_{n2}-g_{k2}\gamma_{n1}\\
&x_{55}=\gamma_{51}\gamma_{52}, &&x_{66}=\gamma_{61}\gamma_{62}&&x_{56}=\gamma_{51}\gamma_{62}+\gamma_{61}\gamma_{52}
\end{align*} with the homomorphism
$$\begin{CD}
\C[\rGr]@>>>\C[\rSL(4|2)]\\
q_{ij}, \lambda_{kn}@>>>y_{ij}, \eta_{kn}\\
a_{55},a_{66},a_{56}@>>>x_{55},x_{66},x_{56}
\end{CD}
$$

\end{proposition}
{\it Proof.} The proof uses an argument similar to the one used to obtain (\ref{wedgeproduct}). Instead of taking the vectors $a$ and $b$ we have to take the first two columns of the matrix (\ref{groupmatrix}).

\hfill$\blacksquare$

 \section{The big cell}\label{bigcell-sec}
A $(2|0)$ subspace of $\C^{4|2}$ is given as the linear span  of two even vectors\footnote{Here we are using implicitly the formalism of the {\it functor of points} to describe a super variety \cite{dm}.}
\be\label{foptsV}
V(A)= \rspan\begin{pmatrix}u_1&v_1\\u_2&v_2\\u_3&v_3\\u_4 &v_4\\\hline\mu_1&\nu_1\\\mu_2&\nu_2\end{pmatrix}\qquad u_i, v_i\in A_0,\; \mu_i,\nu_i\in A_1\,,
\ee
where $A$ is any superalgebra. Clearly there is a right action of $\rGL_2(A)$ over $V(A)$ (change of basis). We assume now that
\be\det\begin{pmatrix}u_1&v_1\\u_2&v_2\end{pmatrix} \quad \hbox{invertible in}
\quad A_0
\label{bigcell}\ee  
{This is a supervariety, which is
an affine open set into the affine superspace $\mathbb{A}^{8|4}$.
It is represented by the superring:
$$
\C[S]=\C[a_{ij},\al_{kl}][T]/((a_{11}a_{22}-a_{12}a_{21})T-1), \,\,
1 \leq i<j \leq 4, \, 5 \leq k<l \leq 6
$$
The condition of invertibility of the determinant function
$a_{11}a_{22}-a_{12}a_{21}$ accounts for the condition in (\ref{bigcell}).}

\medskip
Then, using the right action of $\rGL_2(A)$ we can bring 
{(\ref{foptsV})} to the standard form
$$V(A)=\rspan\begin{pmatrix}\id_{2\times 2}\\P_{2\times 2}\\\psi_{2\times 2}^t\end{pmatrix}\,.
$$
$P$ and $\psi$ are even and odd coordinates in the open subset of $\rGr$ characterized by
(\ref{bigcell}) called the big cell of $\rGr$.  
As for the $N=1$ case, one can show that the subgroup of $\rSL(4|2)$ that leaves invariant the big cell contains the (complexified) $N=2$ super Poincar\'{e} group times the $R$-symmetry (dilations for $N=1$). In fact, we call $\rGr$ the antichiral\footnote{See \cite{fl} to see why this space is the antichiral one.} conformal superspace, while the big cell is the antichiral Minkowski superspace. In this respect we do not coincide with the notation of \cite{hl}-\cite{haho}, where they call directly Minkowski superspace to the Grassmannian.

We do not extend here on this construction, but the condition (\ref{bigcell}) will be also used in the quantum setting.

\section{The quantum Grassmannian} \label{quantumgrassmannian-sec}

We follow  \cite{ma2} to define the quantum group $\rSL_q(r|s)$.

\begin{definition} \label{Manin-def} The  quantum matrix superalgebra $\rM_q(r|s)$ is defined as
$$
\rM_q(r|s)=_{\mathrm{def}}\C_q \langle z_{ij},\xi_{kl}\rangle/\cI_M
$$
where $\C_q\langle z_{ij},\xi_{kl}\rangle$ denotes the free
superalgebra over $\C_q=\C[q,q^{-1}]$
generated by the even variables
$$z_{ij},\qquad  \hbox{ for }\quad 1 \leq i,j \leq r \quad \hbox{ or } \qquad r+1 \leq i,j \leq r+s.$$
and by the odd variables
\begin{align*}&\xi_{kl}&&  \hbox{for  }\quad 1 \leq k \leq r, \quad r+1 \leq l \leq r+s \\&&&\hbox{or   }\,  r+1 \leq k \leq r+s, \quad 1 \leq l \leq r,\end{align*}
satisfying the relations $\xi_{kl}^2=0$ and $\cI_M$ is an ideal that we describe below. We can visualize the generators as a matrix
\be\begin{pmatrix}z_{r\times r}&\xi_{r\times s}\\
\xi_{s\times r}&z_{s\times s}\end{pmatrix}\,.\label{generatorsmatrix}\ee

It is convenient sometimes to have a common notation for even and odd variables.
$$
a_{ij}=\begin{cases} z_{ij} & 1 \leq i,j \leq r, \, \hbox{   or   } \quad
                         r+1 \leq i,j \leq r+s ,\\ \\
              \xi_{ij} &  1 \leq i \leq r,\quad  r+1 \leq j \leq r+s, \, \hbox{   or } \\&
                      r+1 \leq i \leq r+s, \quad 1 \leq j \leq r\,.
\end{cases}
$$

 We  assign a parity to the indices: $p(i)=0$ if $1 \leq i \leq r$ and  $p(i)=1$ if  $r+1 \leq i \leq r+s$. The parity of $a_{ij}$ is  $\pi(a_{ij})=p(i)+p(j)$ mod 2. Then, the ideal $\cI_M$ is generated by the relations \cite{ma2}:

\begin{align}
&a_{ij}a_{il}=(-1)^{\pi(a_{ij})\pi(a_{il})}
q^{(-1)^{p(i)+1}}a_{il}a_{ij}, && \hbox{for  } j < l \nonumber\\&&& \nonumber\\
&a_{ij}a_{kj}=(-1)^{\pi(a_{ij})\pi(a_{kj})}
q^{(-1)^{p(j)+1}}a_{kj}a_{ij}, && \hbox{for  } i < k \nonumber\\ &&\nonumber&\\
&a_{ij}a_{kl}=(-1)^{\pi(a_{ij})\pi(a_{kl})}a_{kl}a_{ij}, &&  \hbox{for  }
i< k,\;j > l \nonumber\\&&&\hbox{or } i > k,\; j < l \nonumber\\&&& \nonumber\\
&a_{ij}a_{kl}-(-1)^{\pi(a_{ij})\pi(a_{kl})}a_{kl}a_{ij}=(-1)^{\pi(a_{ij})\pi(a_{kl})}(q^{-1}-q)
a_{kj}a_{il},&\nonumber&\\
&&& \hbox{for }i<k,\nonumber
j<l\nonumber
\\&&&\label{ManinCR}
\end{align}

$\hbs$
\end{definition}

There is also a comultiplication
$$\begin{CD}\rM_q(m|n)@>\Delta>>\rM_q(m|n)\otimes \rM_q(m|n)\end{CD}$$
$\Delta( a_{ij})=\sum_ka_{ik}\otimes a_{kj}$ and a counit $\varepsilon (a_{ij})=\delta_{ij}$.

One can restrict to $\rSL_q(r|s)$ by setting the quantum Berezinian to 1. The antipode is the usual one (see \cite{ma2} or \cite{fl}, Appendix E). Then $\rSL_q(m|n)$ is a super Hopf algebra.


We can now define the quantum Grassmannian $\rGr_q$ mimicking Proposition \ref{generators-prop}:

\begin{definition}The quantum super Grassmannian $\rGr_q:=Gr_q(2|0,4|2)$ is  the subalgebra of $\rSL_q(4|2)$ generated by the elements\label{qgrassmannian-def}
\begin{align*}
&D_{ij}:=a_{i1}a_{j2}-q^{-1}a_{i2}a_{j1} &&D_{in}:=a_{i1}a_{n2}-q^{-1}a_{i2}a_{n1}&&\\
&D_{55}:=a_{51}a_{52} &&D_{66}:=a_{61}a_{62}\\&D_{56}=a_{51}a_{62}-q^{-1}a_{52}a_{61}\\
\end{align*}
with $1\leq i<j\leq 4$ and $\;n=5,6 $.

\end{definition}

\hfill$\blacksquare$

We want to give a presentation in terms of generators and relations, as in Proposition \ref{superring-prop} for the classical case. Note that, first of all, we have to compute the commutation rules among the $D$'s.  After some (tedious) calculations we arrive at:

\begin{itemize}

\item Let $1\leq i,j,k,l\leq6$ be not all distinct, and $D_{ij}$, $D_{kl}$ not both odd. Then
\be D_{ij}D_{kl}=q^{-1}D_{kl}D_{ij},\qquad (i,j)<(k,l),\; i<j, k<l\,,\label{CR3}\ee where the ordering `$<$' of pairs is the lexicographical ordering.

\item Let $1\leq i,j,k,l\leq6$ be all distinct, and $D_{ij}$, $D_{kl}$ not both odd and $D_{ij},D_{kl}\neq D_{56}$. Then
    \begin{align}
    &D_{ij}D_{kl}=q^{-2}D_{kl}D_{ij}, &&1\leq i<j<k<l\leq 6,\nonumber\\
    &D_{ij}D_{kl}=q^{-2}D_{kl}D_{ij}-(q^{-1}-q)D_{ik}D_{jl} &&1\leq i<k<j<l\leq 6,\nonumber\\
    &D_{ij}D_{kl}=D_{kl}D_{ij}&&1\leq i<k<l<j\leq 6,\label{CR1}
     \end{align}
\item Let  $1\leq i<j\leq4$, $5\leq n\leq m\leq 6$. Then
\begin{align}
&D_{in}D_{jn}= -q^{-1}D_{jn}D_{in}-(q^{-1}-q)D_{ij}D_{nn}=-qD_{jn}D_{in},\nonumber\\
&D_{ij}D_{nm}=q^{-2}D_{nm}D_{ij},\nonumber \\
&D_{i5}D_{j6}=-q^{-2}D_{j6}D_{i5}-(q^{-1}-q)D_{ij}D_{56},\nonumber\\
&D_{i6}D_{j5}=-D_{j5}D_{i6},\nonumber\\
&D_{i5}D_{i6}=-q^{-1}D_{i6}D_{i5},\nonumber\\
&D_{i5}D_{i6}=-q^{-1}D_{i6}D_{i5},
\nonumber\\
&D_{55}D_{66}=q^{-2}D_{66}D_{55}\,,\nonumber \\
&D_{55}D_{56}=0\,. \label{CR2}
\end{align}
\end{itemize}
The Pl\"{u}cker relations are modified. One has for $1\leq i<j<k\leq 4$ and $ n=5, 6$:

\begin{align}
&D_{12}D_{34}-q^{-1}D_{13}D_{24}+q^{-2}D_{14}D_{23}=0, \nonumber\\
&D_{ij}D_{kn}-q^{-1}D_{ik}D_{jn}+q^{-2}D_{jk}D_{in}=0, \nonumber
\\
&D_{i5}D_{j6}+q^{-1}D_{i6}D_{j5}=qD_{ij}D_{56}, \nonumber\\
&D_{in}D_{jn}=qD_{ij}D_{nn}, \nonumber \\
&D_{in}D_{nn}=0, \nonumber\\&D_{i5}D_{66}=-q^{-1}D_{i6}D_{56},\nonumber\\
&D_{i6}D_{55}=-q^{2}D_{i5}D_{56} \nonumber\\
&D_{nn}^2=0,\nonumber\\&D_{55}D_{56}=0, \nonumber \\&D_{66}D_{56}=0\nonumber\\\
&D_{56}D_{56}=(q^{-1}-3q)D_{55}D_{66}\,.\label{quantumplucker}
\end{align}

The first relation in (\ref{CR2}) has been simplified with the use of the fourth relation in (\ref{quantumplucker}).

We have the following

\begin{proposition}The quantum Grassmannian superring $\rGr_q=Gr_q(2|0,4|2)$ is given in terms of generators and relations as
$$\rGr_q=\C_q\langle X_{ij}, X_{mn}, X_{im} \rangle, \cI_q,\qquad 1\leq i<j\leq 4;\quad 5\leq m\leq n\leq 6 \,,$$
where $\cI_q$ is the ideal generated by the commutation relations (\ref{CR3}),(\ref{CR1}),(\ref{CR2}) and the quantum super Pl\"{u}cker relations (\ref{quantumplucker}).

\end{proposition}
{\it Proof.}
We give a sketch of the argument, whose idea is expressed
in \cite{fi1} Theorem 5.4 and also in \cite{fl} Chapter 4.

The super Pl\"{u}cker relations are all the relations  satisfied by the quantum determinants: suppose that there is an extra relation R. Then $R=(q-1)R^{(1)}$.  Then $R^{(1)}$ may be of the form $R^{(1)}=(q-1)R^{(2)}$ or $R^{(1)}$ mod $(q-1)$ not identically 0. In the second case,  since $\rSL_q(m|n)$ is an algebra without torsion, we would have an additional classical Pl\"{u}cker relation, which cannot be. In the first case we have the same possibilities for $R^{(2)}$. At the end of the procedure we will obtain $R^{(n)}=0$, that would be a new classical Pl\"{u}cker condition. But we know that this is not possible.

\hfill$\blacksquare$

\section{The quantum super Grassmannian as a quantum homogeneous superspace}\label{homogenous-sec}
To finish the interpretation of the quantum super Grassmannian as an homogeneous superspace under the quantum supergroup $\rSL_q(4|2)$ we have to see how it is the coaction on $\rGr_q$.
 This is done in the following

 \begin{proposition} The restriction of the comultiplication in $\rSL_q(4|2)$
 $$\begin{CD}
 \rSL_q(4|2)@> \Delta>> \rSL_q(4|2)\otimes\rSL_q(4|2)\\
 a_{ij}@>>> \Delta(a_{ij})=\sum_{k=1}^6a_{ik}\otimes a_{kj}
 \end{CD}
 $$
 to the subalgebra $\rGr_q$ is of the form\footnote{We denote with $\Delta$ both, the comultiplication and its restriction to $\rGr_q$ in order not to burden the notation. The meaning should be clear from the context.}
$$\begin{CD}
 \rGr_q@> \Delta>> \rSL_q(4|2)\otimes\rGr_q\,.
 \end{CD}
 $$

  \end{proposition}
 {\it Proof.} The coaction property is guaranteed by the associativity of the coproduct, so we only have to check that
 $$D_{ij}, D_{im}, D_{mn}\in \rSL_q(4|2)\otimes\rGr_q\,.$$

 Let us denote as $D_{ij}^{kl}=a_{ik}a_{jl}-q^{-1}a_{il}a_{jk}$, so in the previous notation $D_{ij}=D_{ij}^{12}$. After some calculations one can prove

 \begin{enumerate}
 \item Let us call $P$ the condition $1\leq k,l\leq 6$ and  at least one of the two indices is less that 5. For $1\leq i<j \leq 4$:

\begin{align*}\Delta(D_{ij})&= \sum_{P \cap (k<l)} D_{ij}^{kl} \otimes D_{kl}^{12}- (a_{i5}a_{j6}+q^{-1}a_{i6}a_{j5})  \otimes D_{56}\\&- (1+q^{-2})\sum_{5\leq k \leq 6}a_{ik}a_{jk} \otimes D_{kk}\,.\end{align*}

\item For $1\leq i \leq 4$ and $5\leq m \leq 6$:

\begin{align*}
\Delta(D_{im})&=\sum_{\substack{k<5 \\ k<l}} a_{ik}a_{ml} \otimes D_{kl} -q^{-1}\sum_{\substack{k<5 \\ l<k}} a_{ik}a_{ml} \otimes D_{lk}\\& + (a_{i5}a_{m6}+q^{-1}a_{i6}a_{m5}) \otimes D_{56}\\
& + (1+q^{-2})\sum_{5\leq k \leq 6} a_{ik}a_{mk} \otimes D_{kk}+q^{-1} \sum_{\substack{k\geq5 \\ l<5}} a_{ik}a_{ml} \otimes D_{lk}\,.
\end{align*}

\item  For $5 \leq m,n \leq 6$:
\begin{align*}
\Delta(D_{56})&=\sum_{\substack{k<5 \\ k<l}} a_{5k}a_{6l} \otimes D_{kl} -q^{-1}\sum_{\substack{k<5 \\ l<k}} a_{5k}a_{6l} \otimes D_{lk}\\& + (a_{55}a_{66}+q^{-1}a_{56}a_{65}) \otimes D_{56}\\
& + (1+q^{-2})\sum_{5\leq k \leq 6} a_{5k}a_{6k} \otimes D_{kk}+q^{-1} \sum_{\substack{k\geq5 \\ l<5}} a_{5k}a_{6l} \otimes D_{lk}\,,
\end{align*}

and
\begin{align*}
\Delta(D_{nn})= \sum_{1\leq k<l \leq 6} a_{nk}a_{nl} \otimes D_{kl} +\sum_{5\leq k \leq 6} a_{nk}^{2} \otimes D_{kk}\,.
\end{align*}
\end{enumerate}

This proves our statement.

\hfill$\blacksquare$


\section{Quantum super bundles: quantum chiral Minkowski
superspace}\label{quantumbundle-sec}
In this section we want to reinterpret our construction in
the framework of quantum principal bundles, as in \cite{afl} and references therein.
We shall concentrate our attention on the local picture, that is,
we want to look at the quantization of a super
bundle $S \lra \C^{4|4}$, with base space 
the chiral Minkowski superspace $\C^{4|4}$, which we interpret as the
big cell into the Grassmannian supermanifold $\mathrm{Gr}$
{(see also Sec. \ref{bigcell-sec}).}

We shall not develop a full theory of quantum principal super bundles,
but we will recall the key definitions in order to put in the correct
framework our construction.

We start with the classical definition.

\begin{definition}\label{p-princ}
Let $X$ and $M$ be topological spaces, $P$ a topological group and
$\wp: X \lra M$ a continuous function. We say
that $(X, M, \wp, P)$ is a \textit{$P$-principal bundle}
(or {\it principal bundle} for short) with total space
$X$ and base $M$, if the following conditions hold
\begin{enumerate}
\item $\wp$ is surjective.
\item $P$ acts freely from the right on $X$.
\item $P$ acts transitively on the fiber $\wp^{-1}(m)$
of each point $m \in M$.
\item $X$ is locally trivial over $M$, i.e. there is an open covering
$M=\cup U_i$ and homeomorphisms
$\sigma_i:\wp^{-1}(U_i)  \lra  U_i \times P$
that are $P$-equivariant i.e.,
$\sigma_i(up)=\sigma_i(u)p$, $u \in U_i$, $p \in P$.
\end{enumerate}

\medskip

If $X\cong M\times P$ we say that the bundle is globally trivial.

\hfill$\blacksquare$
\end{definition}

We can then define {algebraic, analytic or smooth} $P$-principal
bundles, by the taking objects and morphisms in the appropriate
categories. There is clearly no obstacle in writing the same
definition in the super context, provided we exert some care
in the definition of surjectivity
(see \cite{ccf}, Section 8.1 for details).
We would like, however, to take a different
route.

We turn to the notion of Hopf-Galois extension, that is most fruitful
for the quantization. Our definition in the super category is
the same as for the ordinary one (see \cite{mont} for more details in the ordinary category).

\begin{definition} \label{hg-def}
Let $(H, \Delta, \epsilon, S)$ be a Hopf superalgebra and
$A$ be an $H$-comodule superalgebra with coaction
$\delta : A \lra A \otimes H$.
Let
\beq \label{coinvariant-def}
B:=A^{\coinv\,H} := \{a \in A \,|\, \delta(a) = a \otimes 1\}~.
\eeq
The  extension $A$ of the superalgebra $B$ is called
$H$-{Hopf-Galois} (or simply {Hopf-Galois}) if the
map
$$
\chi:A \otimes_B A \lra A \otimes H, \qquad \chi=(m_A \otimes
\mathrm{id})(\mathrm{id}
\otimes_B \delta)
$$
called the {canonical map}, is bijective ($m_A$ denotes the
multiplication in $A$).

The extension $B=A^{\coinv\,H}\subset A$ is called
$H$-\textit{principal comodule superalgebra} if it is Hopf-Galois and
$A$ is $H$-equivariantly projective as a left $B$-supermodule, i.e., there
exists a left $B$-supermodule and right $H$-comodule morphism $s : A \to B
\otimes A$ that is a section of the (restricted) product $m : B
\otimes A \to A$.

In particular if $H$ is a Hopf algebra with bijective antipode
over a field, the condition of equivariant projectivity of $A$ is
equivalent to that of faithful flatness of $A$ (see \cite{ss},
\cite{abps}).
\hfill$\blacksquare$
\end{definition}

{We now follow \cite{afl} Sec. 2, in giving the definition  of quantum
principal bundle, though it differs slightly from the one given
in the literature, which also requires the existence of a differential
structure (see e.g. \cite{bm} Ch. 5). We plan to explore such
structures in a forthcoming paper.

\begin{definition}\label{pb-def}
We define \textit{quantum principal bundle} 
a pair $(A,B)$, where $A$ is an $H$-Hopf Galois extension and
$A$ is $H$-equivariantly projective as a left $B$-supermodule, that is,
$A$ is an $H$-principal comodule superalgebra.

\hfill$\blacksquare$
\end{definition}}

In the ordinary case, the notion of affine $H$-principal bundle
coincides with Definition \ref{p-princ} when we take $H=\cO(P)$,
$A=\cO(X)$ and $B=\cO(M)$, where $\cO(X)$ denotes the algebra
of functions on $X$ (algebraic, differential, holomorphic, etc).
The Hopf-Galois condition is equivalent to saying that the action of $P$
on $X$ is free and the equivariance property means that the bundle
is locally trivial.

We assume, in the algebraic setting, that all our varieties are affine.

\medskip
There is a special case of Hopf-Galois extensions, corresponding to
a globally trivial principal bundle. In this case the
technical conditions of Definition \ref{pb-def} are
  automatically satisfied. We shall focus on this case
leaving aside the general one.

\begin{definition}\label{cleft-def}
Let $H$ be a Hopf superalgebra and $A$ an $H$-comodule superalgebra.
The algebra extension $A^{\coinv\,H}\subset A$ is called
a \textit{cleft extension} if
there is a right $H$-comodule map $j : H \to A$, called \textit{cleaving map},
that is convolution invertible, i.e. there exists a map $h:H\to A$  
such that the convolution product $j\star h$ satisfies,
$$j\star h:=m_A\circ(j\otimes h)\circ\Delta(f)=\epsilon(f)\cdot 1\, $$
or, in  Sweedler notation
$\Delta(f)=\sum f_1 \otimes f_2$,
 $$\sum j(f_1)h(f_2)=\epsilon(f)\cdot 1\,$$
for all $f \in H$. The map $h$ is the convolution inverse of $j$.


An extension $A^{\coinv\,H}\subset A$ is called a
{trivial extension} if there exists such map.

Notice that when $j$ is an algebra map,
its convolution inverse is just $h=j \circ S^{-1}$.

\hfill $\blacksquare$
\end{definition}

A trivial extension is actually a Hopf-Galois extension and a
principal bundle. {When $j$ is an algebra map, 
we have an algebra isomorphism $A \cong B \# H$} (see
\cite{mont}, Sections 4.1 and 7.2 for the smash product `$\,\#\,$'),
which in the classical case means that we have a trivial bundle
(see \cite{mont} Chapter 8 and \cite{bm} Sec. 5.1.2).

\medskip
We now examine an example with physical significance coming
from our previous treatment. Consider the set
of $4 \times 2 \, | \, 2 \times 2$
supermatrices with complex entries

\beq
\begin{pmatrix} a_{11} &  a_{12} \\a_{21} &a_{22}\\ a_{31} &a_{32}\\a_{41}&a_{41}\\
\alpha_{51}&\alpha_{52}\\\alpha_{61}&\alpha_{62}
\end{pmatrix}\,.
\eeq

 This can be seen as the affine superspace
$\bA^{8|4}$ described by the coordinate superalgebra
$\C[a_{ij},\al_{kl}]$ with $ i=1,\dots, 4$, $\,j,l=1,2$, $\,k=5,6$.
 As in the ordinary
setting, we can view elements in $\bA^{8|4}$ as $2|0$ subspaces
of $\C^{4|2}$:
$$
W=\Span\{a_1,a_2\} \subset \C^{4|2} 
$$
In this way,
$W$ may also be viewed as an element in $\mathrm{Gr}$.

\medskip
In the superspace $\bA^{8|4}$ consider the open subset $S$ consisting
of matrices such that the minor formed with $a_{ij}$, $i,j=1,2$ is invertible.
This open subset $S$ is described by its coordinate superalgebra:
$$
\C[S]=\C[a_{ij},\al_{kl}][T]/((a_{11}a_{22}-a_{12}a_{21})T-1)
$$
We have a right action of $\rGL_2(\C)$ on $S$ corresponding to
the change of basis of such subspaces:
$$
\Span\{a_1,a_2\}, g \mapsto \Span\{a_1 \cdot g, a_2 \cdot g\},\qquad g\in \rGL_2(\C)\,.
$$

\begin{proposition}\label{affine-prop}
Let the notation be as above.
The quotient of $S$
by the right $\rGL_2(\C)$ action is
an affine superspace $\bM$ of dimension $4|4$, the $N=2$ chiral Minkowski superspace $\bM$.
\end{proposition}

{\it Proof.}
We can write:
\beq
\bM=\left\{ (a_1,a_2),\, a_1,a_2 \in \C^{4|2} \;|\; \det\begin{pmatrix}  a_{11} &  a_{12} \\ a_{21} &  a_{22}
\end{pmatrix} \neq 0 \right\}
\Big/\rGL_2(\C) \qquad
\eeq
In the quotient $\bM$ we can choose a (unique) representative
$(u,v)$ for $(a_1,a_2)$ of the form:
\beq
\left\{
\begin{pmatrix} 1 \\  0 \\u_{1} \\u_{2}\\ \nu_{3} \\ \nu_{4}\end{pmatrix},
\begin{pmatrix} 0 \\  1 \\v_{1} \\v_{2} \\  \eta_{3} \\ \eta_{4}
\end{pmatrix}\right\}\,,
\eeq
so $\bM$ is $\C^{4|4}$.

\hfill$\blacksquare$

We notice that $\bM$ is naturally identified with the dense open
set of the Grassmannian
$\mathrm{Gr}$ in the Pl\"ucker embedding, determined by the
invertibility of the coordinate $q_{12}$ in $\mathbb{P}^{8|8}$.%

\medskip
We now would like to retrieve a set of global coordinates for
$\bM$ starting from the global coordinates $a_{ij}$  for $S$.
Let $\C[\rGL_2]=\C[g_{ij}][T]/((g_{11}g_{22}-g_{12}g_{21})T-1)$ be the
coordinate algebra for the algebraic group $\rGL_2(\C)$.
Let us write heuristically the equation relating the generators
of $\C[S]$, $\C[\rGL_2]$ and  the polynomial
superalgebra $\C[\bM]:=\C[u_{ij}, \nu_{kl}]$
\beq\label{mink-coord}
\begin{pmatrix}
a_{11} & a_{12} \\
a_{21} & a_{22} \\
a_{31} & a_{32} \\
a_{41} & a_{42} \\
\al_{51} & \al_{52} \\
\al_{61} & \al_{62}
\end{pmatrix}
=
\begin{pmatrix} 1 &  0 \\0 &  1 \\
u_{31} & u_{32}\\
u_{41} & u_{42}\\
\nu_{51}& \nu_{52} \\
\nu_{61} & \nu_{62} \end{pmatrix}
\begin{pmatrix}  g_{11} &  g_{12} \\ g_{21} &  g_{22}
\end{pmatrix}
\eeq
We obtain immediately:
$$
\begin{pmatrix}  g_{11} &  g_{12} \\ g_{21} &  g_{22}
\end{pmatrix} = \begin{pmatrix}  a_{11} &  a_{12} \\ a_{21} &  a_{22}
\end{pmatrix}
$$
and then with a short calculation,
$$
u_{i1}= -d_{2i}d_{12}^{-1}\qquad u_{i2}= d_{1i}d_{12}^{-1}$$
$$\nu_{k1}= -d_{2k}d_{12}^{-1}\qquad \nu_{k2}= d_{1k}d_{12}^{-1}
$$
for $i=3,4$ and $k=5,6$, where
$$d_{rs}:=a_{r1}a_{s2}-a_{r2}a_{s1} \qquad r<s\,.$$

\begin{proposition}Let the notation be as above.

 \begin{enumerate}
\item The complex supermanifold $S$ is diffeomorphic to the supermanifold
$\C^{4|4} \times \rGL_2(\C)$:

$$\begin{CD} S @>\psi>> \C^{4|4} \times \rGL_2(\C)\,,\end{CD}$$
with
\begin{align*}
& \psi^*(g_{ij})=a_{ij} \\
&\psi^*(u_{i1})= -d_{2i}d_{12}^{-1}, &&
\psi^*(u_{i2})=  d_{1i}d_{12}^{-1}.\\
&\psi^*(\nu_{k1})= -d_{2k}d_{12}^{-1} &&
\psi^*(\nu_{k2})= d_{1k}d_{12}^{-1}
\end{align*}

\item The diffeomorphism $\psi$ is $\rGL_2(\C)$-equivariant
with respect to the right $\rGL_2(\C)$ action, hence
$S/\rGL_2(\C) \cong \C^{4|4}$.

\end{enumerate}
\end{proposition}

{\it Proof.} We notice that $\psi$ is invertible,
$\psi^{-1}$ is given by:
$$
(\psi^{-1})^*(a_{ij})=g_{ij}
$$
and the rest follows from equation (\ref{mink-coord}).

The right equivariance of $\psi$ is a simple calculation, taking into account that the determinants $d_{ij}$ transform as $d_{ij}\det g'$, were $g'\in\rGL_2(\C)$.

\hfill$\blacksquare$

\begin{lemma}\label{coinvariant-lemma}
The coordinate superalgebra $\C[\bM]:=\C[u_{ij}, \nu_{kl}]$
is isomorphic to $\C[S]^{\coinv \,\C[\rGL_2]}$ the coinvariants in
$\C[S]$ with respect to the $\C[\rGL_2]$ right coaction $\de$:
$$
\begin{CD} \C[S] @>\de>> \C[S] \otimes \C[\rGL_2] \\
\begin{pmatrix}
a_{11} & a_{12} \\
a_{21} & a_{22} \\
a_{31} & a_{32} \\
a_{41} & a_{42} \\
\al_{51} & \al_{52} \\
\al_{61} & \al_{62}
\end{pmatrix} @>>>\begin{pmatrix}
a_{11} & a_{12} \\
a_{21} & a_{22} \\
a_{31} & a_{32} \\
a_{41} & a_{42} \\
\al_{51} & \al_{52} \\
\al_{61} & \al_{62}
\end{pmatrix} \otimes \begin{pmatrix}  g_{11} &  g_{12} \\ g_{21} &  g_{22}
\end{pmatrix}
\end{CD}
$$
\end{lemma}

{\it Proof}.
In our heuristic calculation we computed expressions for
the coordinates on $\bM$. We claim that
these are coinvariants, so we need to show $\de(c)=c\otimes 1$
for any $c\in\{u_{ij}, \nu_{kl}\}$. A little calculation gives us
$$\de(d_{rs})= d_{rs} \otimes (g_{11}g_{22}-g_{12}g_{21})\; \Rightarrow\;\de(d_{rs}d_{12}^{-1})= d_{rs}d_{12}^{-1} \otimes 1$$
which proves our claim.

On the other hand,  the space of functions on $S$ that are invariant under the right translation of the group can be identified with  the space of functions on the quotient $S/\rGL_2(\C)$. Since we have global coordinates in $\bM$, any other invariant will be a function of these coordinates. In the Hopf algebra language, this means that $\{u_{ij}, \nu_{kl}\}$ are the only independent coinvariants.

\hfill$\blacksquare$

\begin{proposition}\label{cleaving-prop}
Let the notation be as above.
The natural projection $p:S \lra S/\rGL_2(\C)$ is a trivial principal
bundle on the chiral superspace.
\end{proposition}

{\it Proof.} We will show that $\C[S]$ is a trivial $\C[\rGL_2]$-Hopf
Galois extension of $\C[\bM]$.
In our previous lemma we proved that $\C[\bM]\cong
\C[S]^{\coinv \,\C[\rGL_2]}$, so if we give an algebra
cleaving map we will have proven our proposition. We define

$$\begin{CD}
 \C[\rGL_2] @>j>> \C[S]\\
g_{ij}@>>>a_{ij}\,.\end{CD}
$$
We leave to the reader the easy check that $j$ is convolution invertible with convolution inverse
$$h=j\circ S\,.$$ Moreover, the calculation below shows that $j$ is a $\C[\rGL_2(\C)]$-comodule map,
$$(\delta \circ j) (g_{ij})=\delta(a_{ij})=\sum a_{ik}\otimes g_{kj}.$$
$$((j \otimes \rid) \circ \Delta) (g_{ij})=(j \otimes \rid)(\sum g_{ik} \otimes g_{kj})=\sum a_{ik}\otimes g_{kj}.$$

$$\Rightarrow~~~\delta \circ j= (j \otimes \rid) \circ \Delta.$$
This proves the result.

\hfill$\blacksquare$

We now go to the quantum setting, where
we lose the geometric
interpretation and we retain only the algebraic point of view.
Hence a quantum principal super bundle over an affine base
is just understood as a Hopf-Galois extension with the properties mentioned in Definition \ref{pb-def}.

\medskip

We want to study the quantization of the example studied above.
Let $\C_q[S]$ be the quantization of the superalgebra $\C[S]$
obtained by taking the Manin relations (\ref{ManinCR}) among the entries still
denoted as $a_{ij}$ and $\al_{kl}$, with $i,j=1,\dots,4$ and $k,l=5,6$.

\begin{definition} The $N=2$ quantum chiral Mikowski superspace,
 $\C_q[\bM]$, is  the superalgebra generated by the elements

$$
\Tilde{u}_{i1}= -q^{-1}D_{2i}D_{12}^{-1},\qquad \Tilde{u}_{i2}= D_{1i}D_{12}^{-1},$$
$$\Tilde{\nu}_{k1}= -q^{-1}D_{2k}D_{12}^{-1},\qquad \Tilde{\nu}_{k2}= D_{1k}D_{12}^{-1}\,,
$$
for $i=3,4$ and $k=5,6$ in $\C_q[\rGr]$, where
$$
D_{rs}:= a_{r1}a_{s2}-q^{-1}a_{r2}a_{s1}, \qquad r<s.
$$

\end{definition}

\hfill$\blacksquare$

Using our previous computations for commutation relations among $D_{rs}$'s we get the following commutation relations among $\Tilde{u}_{ij}$'s and $\Tilde{\nu}_{kl}$'s:

\begin{align}
&\Tilde{u}_{i2}\Tilde{u}_{i1}= q^{-1}\Tilde{u}_{i1}\Tilde{u}_{i2}, \qquad &&i=3,4, \nonumber \\
&\Tilde{\nu}_{k1}\Tilde{\nu}_{k2}=-q^{-1}\Tilde{\nu}_{k2}\Tilde{\nu}_{k1}, \qquad &&k=5,6, \nonumber \\
&\Tilde{\nu}_{5l}\Tilde{\nu}_{6l}=-q^{-1}\Tilde{\nu}_{6l}\Tilde{\nu}_{5l}, \qquad &&l=1,2,\nonumber \\
&\Tilde{u}_{3j}\Tilde{u}_{4j}=q^{-1}\Tilde{u}_{4j}\Tilde{u}_{3j}, \qquad &&j=1,2,\nonumber \\
&\Tilde{u}_{ij}\Tilde{\nu}_{kj}=q^{-1}\Tilde{\nu}_{kj}\Tilde{u}_{ij}, \qquad &&j=1,2~~i=3,4~~k=5,6,\nonumber \\
&\Tilde{u}_{i1}\Tilde{\nu}_{k2}=\Tilde{\nu}_{k2}\Tilde{u}_{i1}, \qquad &&i=3,4~~k=5,6,\nonumber\\
&\Tilde{u}_{31}\Tilde{u}_{42}=\Tilde{u}_{42}\Tilde{u}_{31},\nonumber\\
&\Tilde{\nu}_{51}\Tilde{\nu}_{62}=-\Tilde{\nu}_{62}\Tilde{\nu}_{51}, \nonumber \\
&\Tilde{u}_{32}\Tilde{u}_{41}-\Tilde{u}_{41}\Tilde{u}_{32}= (q^{-1}-q)\Tilde{u}_{42}\Tilde{u}_{31}, \nonumber \\
&\Tilde{u}_{i2}\Tilde{\nu}_{k1}-\Tilde{\nu}_{k1}\Tilde{u}_{i2}= (q^{-1}-q)\Tilde{\nu}_{k2}\Tilde{u}_{i1},\qquad &&i=3,4~~k=5,6,\nonumber\\
&\Tilde{\nu}_{52}\Tilde{\nu}_{61}+\Tilde{\nu}_{61}\Tilde{\nu}_{52}=-(q^{-1}-q)\Tilde{\nu}_{62}\Tilde{\nu}_{51}\,.\label{CRMinkowski}
\end{align}

We have the following

\begin{proposition}
The quantum chiral Minkowski superspace $\C_{q}[M]$ is isomorphic to the quantum superalgebra of matrices $\rM_q(2|2)$ (Definition \ref{Manin-def}).

\end{proposition}

{\it Proof.}
We define the  map $\beta:\C_{q}[M] \longrightarrow \rM_q(2|2)$  by giving it on the generators as follows:
\begin{eqnarray*}
\beta(\tilde{u}_{ij}):=z_{rs}~~~~~~\hbox{where}~~r=i-2~~\hbox{and}~~s=\begin{cases} 1~~\hbox{if}~~j=2,\\
2~~\hbox{if}~~j=1,
\end{cases}
\end{eqnarray*}
\begin{eqnarray*}
\beta(\tilde{\nu}_{kl}):= \xi_{mn}~~~~~~\hbox{where}~~m=k-2~~\hbox{and}~~n=\begin{cases} 1~~\hbox{if}~~l=2,\\
2~~\hbox{if}~~l=1\,.
\end{cases}
\end{eqnarray*}
It is clearly bijective. Comparing  the commutation relations (\ref{ManinCR}) with (\ref{CRMinkowski}) it follows that $\beta$ is an isomorphism.

 \hfill$\blacksquare$

We now want to show the main result for this section.

\begin{theorem}
The quantum superalgebra $\C_q[S]$ is a trivial quantum
principal super bundle on the quantum chiral Minkowski superspace.
\end{theorem}

{\it Proof.}
We need to show that $\C_q[S]$ is a trivial Hopf-Galois extension of $\C_q[M]$. We will proceed similarly to the classical case. It is easy to see that the quantum version of Lemma \ref{coinvariant-lemma} also holds. It is enough to check that
$$\de_{q}(D_{rs})= D_{rs} \otimes (g_{11}g_{22}-q^{-1}g_{12}g_{21})\,.$$
Therefore, we need to give an algebra cleaving map $j_{q}: \C_{q}[\rGL_2(\C)] \lra \C_q[S]$. Define:
$$
j_{q}(g_{ij}):=a_{ij},\qquad h_q=j_q\circ S_q,
$$
so $h_{q}:\C_{q}[\rGL_2(\C)] \lra \C_{q}[S]$,
$$h_{q}(g_{11}):=D^{-1}a_{22} \qquad  h_{q}(g_{12}):=-qD^{-1}a_{12}$$
$$h_{q}(g_{21}):=-q^{-1}D^{-1}a_{21} \qquad h_{q}(g_{22}):= D^{-1}a_{11},$$
where $D:=a_{11}a_{22}-q^{-1}a_{12}a_{21}$. One can observe that:
$$j_{q} \star h_{q}= \varepsilon.1= h_{q}\star j_{q},$$
where $\star$ denotes the convolution product, i.e $j_{q}$ is convolution invertible. Moreover, a similar calculation to the one given in Proposition \ref{cleaving-prop} shows that $j_{q}$ is a $\C_{q}[\rGL_{2}]$-comodule map, i.e. $\delta_{q} \circ j_{q}= (j_{q} \otimes\rid) \circ \Delta$. Therefore, $j_{q}$ is an algebra cleaving map and $\C_{q}[M]\subset \C_{q}[S]$ is a trivial extension. Hence the theorem is proven.\hfill$\blacksquare$

\section*{Acknowledgements}

We thank P. Aschieri, E. Latini, T. Weber, R. Buachalla, A. Carotenuto, E. Shemyakova and T. Voronov for helpful observations. 

This work is supported by the Spanish Grants FIS2017-84440-C2- 1-P funded by MCIN/AEI/10.13039/501100011033 “ERDF A way of making Europe”, PID2020-116567GB-C21 funded by MCIN/AEI/10.13039/501100011033 and the project PROMETEO/2020/079 (Generalitat Valenciana).

\hfill\eject

\end{document}